\documentclass[15pt, a4paper]{article}
\usepackage{}
\usepackage{amsfonts}
\usepackage{mathbbold}
\usepackage{bbm}
\usepackage{mathrsfs}
\usepackage{latexsym}
\usepackage{graphicx}
\usepackage{amssymb}
\usepackage{cite}

\newtheorem{theorem}{Theorem}[section]
\newtheorem{lemma}[theorem]{Lemma}
\newtheorem{corollary}[theorem]{Corollary}

\def\R2{\par\noindent{\bf Remark 2~}}

\title{{\Large \bf On the maximum spectral radius of planar graphs\thanks{Supported by National natural science foundation of China (NSFC)
(Nos. Nos. 12171222, 12101285), Natural science foundation of Guangdong province (No. 2021A1515010254).}}}
\author{Guanglong Yu$^{a}$
~ Lin Sun$^{b}$\thanks{Corresponding authors, E-mail addresses:
yglong01@163.com (G. Yu).} ~
\\ ~ \\
{\footnotesize $^a$School of Mathematics and Systems Science, Guangdong Polytechnic Normal University,}\\ {\footnotesize  Guangzhou, Guangdong, 510665, P.R. China}\\
{\footnotesize $^b$Department of Mathematics, Lingnan Normal
University,  Zhanjiang, Guangdong, 524048, P.R. China}}
\date{}

\begin{document}
\maketitle

\begin{abstract}
This paper investigates the maximum spectral radius of planar graphs with concrete fixed number of vertices, providing some tight bounds on the maximum spectral radius of general planar graph resorting to its order, and confirming that among all planar graphs containing dominating vertex with concrete fixed order $n \geq 48$, the join of $P_{2}$ and $P_{n-2}$ attains the maximum spectral radius.


\end{abstract}
\baselineskip 18.6pt

\section{Introduction}

\ \ \ \ The study of planar graphs has a long history because of its good structural properties, topological properties, algebraic properties, and so on.
To questions in spectral extremal graph theory asking to maximize or minimize eigenvalues
over a fixed family of graphs, Boots and Royle, and independently, Cao and Vince conjectured that
the join of $P_{2}$ and $P_{n-2}$ attains the maximum spectral radius among all planar graphs on $n \geq 9$ vertices \cite{BBGR, DCAV}. This conjecture has been open for more than 30 years although many researchers tries to solve it. Just up to 2017, the study on this conjecture got a substantial development that Tait and Tobin proved the conjecture holding for graphs of sufficiently large order \cite{MTJT}. Unfortunately, the conjecture for concrete fixed order is still open. This paper investigates the maximum spectral radius of planar graphs with concrete fixed order, providing some tight bounds on the maximum spectral radius of general planar graph resorting to its order, and confirming that among all planar graphs containing dominating vertex with fixed order $n \geq 48$, the join of $P_{2}$ and $P_{n-2}$ attains the maximum spectral radius.

\begin{Nont}
\end{Nont}

All graphs considered in this paper are undirected and simple, i.e. neither loop nor multiple edge is allowed. For a set $S$, $|S|$ is employed to denote its cardinality.
And for a $graph$ $G=(V, E)$ consisting of a nonempty vertex set $V=V(G)$ and an edge set $E=E(G)$, $n=|V(G)|$ is called the $order$, $m=|E(G)|$ is called the $edge$ $number$.

Recall that in a graph $G$, for a walk $W=v_{0}e_{0}v_{1}e_{1}\cdots e_{k-1}v_{k}$,
where $e_{i}=v_{i}
v_{i+1}$ for $0\le i\le k-1$ (possibly, $e_{i}=e_{j}$ for $0\leq i\neq j\leq k-1$, $v_{i}=v_{j}$ for $0\leq i\neq j\leq k$), written as
$W=v_{0}v_{1}\cdots v_{k}$
or $W=e_{1}e_{2}\cdots e_{k}$ for short, the integer $k$ is called the $length$.
A $path$ is a walk in which the vertices are pairwise different; a $circuit$ is a closed walk; a $cycle$ is a circuit in which the vertices are pairwise different. A cycle with length $k$ is called a $k$-$cycle$; $k$-$walk$ and $k$-$path$ are defined similarly. We sometime denote by $P_{k}$ ($C_{k}$) a $k$-path (a $k$-cycle) for convenience. A cycle in a graph $G$ is called $Hamilton$ $cycle$ if it contains all vertices of $G$. Denoted by $K_{n}$ a $complete$ $graph$ of order $n$, $K_{s, t}$ a $complete$ $bipartite$ $graph$ with one part of $s$ vertices, the other one of $t$ vertices.

For two graphs $H$ and $G$, $H\cup G$ is employed to denote the the graph obtained by $V(H)\cup V(G)$ and $E(H)\cup E(G)$. In a graph $G$, for two subgraphs $G_{1}$, $G_{2}$ satisfying  $V(G_{1})\cap V(G_{2})=\emptyset$, $E(G_{1}, G_{2})$ is employed to denote the edge set between $V(G_{1})$ and $V(G_{1})$. A vertex $v$ is said to be $incident$ with graph $H$ if $E(H, v)\neq \emptyset$.

Two graphs $H$, $G$ are called $disjoint$ or $independent$ if $V(H)\cap V(G)=\emptyset$, and there is no edge between $H$ and $G$. Denote by $H\bigtriangledown G$ the $join$ of two disjoint graphs $H$ and $G$ obtained from $H\cup G$ by adding edges between every pair of vertices $u$, $v$ where $u\in V(H)$, $v\in V(G)$.

For a graph $G$, we denote by $G+e$ ($G-e$)
the graph obtained from $G$ by adding a new edge $e$ if $e\notin
 E(G)$ (deleting an edge $e\in
 E(G)$); $G[S]$ the subgraph induced by $S\subseteq
V(G)$; $G-S$ the graph obtained from $G$ by deleting all the vertices in $S$ and all the edges incident with the vertices in $S$; $e(S)$ the number of edges in $G[S]$; $G-v$ for $G-\{v\}$ for short. For two subsets $S_{1}\subseteq
V(G), S_{2}\subseteq
V(G)$, denote by $E(S_{1}, S_{2})$ the set of edges between $S_{1}$ and $S_{2}$, $e(S_{1}, S_{2})=|E(S_{1}, S_{2})|$, $G(S_{1}, S_{2})$ the bipartite subgraph with parts $S_{1}, S_{2}$ and edge set $E(S_{1}, S_{2})$. In a graph $G$, a vertex $v$ (an edge $e$) is called a $cut$ $vertex$ ($cut$ $edge$) if the component number of $G-v$ ($G-e$) increases. A 2-connected graph is a connected graph having no cut vertex.

In a graph, we denote by $u\sim v$ for vertex $u$ being adjacent to vertex $v$. For a graph $G$, denoted by $N_{G}(v)$ ($d_{G}(v)=|N_{G}(v)|$) the $neighbor$ $set$ (degree) of vertex $v$; $N_{G}[v]=N_{G}(v)\cup\{v\}$, $G^{v}=G[N_{G}[v]]$; $\Delta(G)$ or $\Delta$ for short ($\delta(G)$ or $\delta$ for short) the $maximal$ $degree$ ($minimal$ $degree$). In a graph $G$ with $n$ vertices, a vertex $v$ is called $dominating$ $vertex$ if $d_{G}(v)=n-1$.

A graph is said to be $planar$ (or embeddable in the plane), if it can be drawn in
the plane so that its edges intersect only at their ends. Such a drawing is called
a $planar$ $embedding$ of the graph (see Fig. 1.1 for example). A planar embedding of a planar graph $G$ partitions the plane into a number of edgewise-connected
open sets. These sets are called the $faces$ of $G$. The number of faces in a planar embedding of a planar graph $G$, denoted by $\mathbbm{f}$ or $\mathbbm{f}_{G}$. Among the faces of a planar embedding, the outer one is called the $outer$ $face$, and any one of the other faces is called the $inner$ $face$ (see Fig. 1.1 for example). It can be seen that the boundary of a face $f$ in a planar embedding of a planar graph, dented by $B(f)$, is a circuit. For a planar graph $G$, in its a planar embedding $\widetilde{G}$, we denote by $O_{\widetilde{G}}$ the outer face. As shown in Fig. 1.1, we can see that $f_{1}-f_{5}$ are inner faces, $O_{\widetilde{G}}=f_{6}$, $B(f_{1})=v_{1}v_{4}v_{5}v_{1}$.

\setlength{\unitlength}{0.7pt}
\begin{center}
\begin{picture}(526,153)
\put(138,81){\circle*{4}}
\put(222,32){\circle*{4}}
\qbezier(138,81)(180,57)(222,32)
\put(311,81){\circle*{4}}
\qbezier(222,32)(266,57)(311,81)
\put(223,133){\circle*{4}}
\qbezier(138,81)(180,107)(223,133)
\qbezier(223,133)(267,107)(311,81)
\qbezier(223,133)(626,77)(222,32)
\qbezier(138,81)(224,81)(311,81)
\put(21,81){\circle*{4}}
\qbezier(223,133)(122,107)(21,81)
\qbezier(21,81)(121,57)(222,32)
\qbezier(138,81)(79,81)(21,81)
\put(219,140){$v_{1}$}
\put(216,20){$v_{2}$}
\put(3,79){$v_{3}$}
\put(130,69){$v_{4}$}
\put(316,79){$v_{5}$}
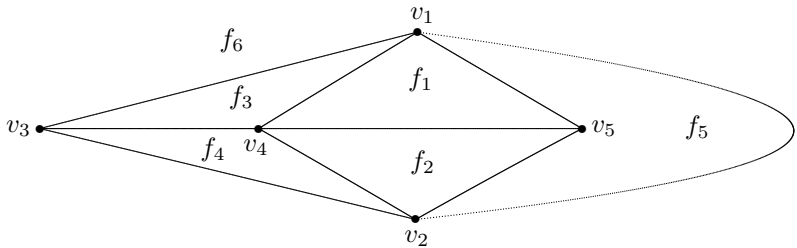
\put(31,-9){Fig. 1.1. a planar graph $G$ and its a planar embedding $\widetilde{G}$}
\put(218,104){$f_{1}$}
\put(219,59){$f_{2}$}
\put(122,94){$f_{3}$}
\put(107,66){$f_{4}$}
\put(366,78){$f_{5}$}
\put(117,125){$f_{6}$}

\end{picture}
\end{center}

A simple planar graph is (edge) $maximal$ if no edge can be added to the graph without violating its simplicity, or planarity.

A graph $H$ is called a $minor$ or $H$-$minor$ of $G$, or $G$ is called
a $H$-$minor$ graph if $H$ can be obtained from $G$ by deleting edges,
contracting edges, and deleting isolated (degree zero) vertices.
Given a graph $H$, a graph $G$ is $H$-$minor$ $free$ (or $H$ $minor$ $free$) if $H$ is not a
minor of $G$.

\begin{Our}
\end{Our}

Denote by $A_{G}$ the $adjacency$ matrix of a graph $G$. It is known that $A_{G}$ is symmetric. The $spectral$ $radius$ (or $A$-$spectral$ $radius$) of graph $G$, denoted by $\rho(G)$, is defined to be the maximum eigenvalue of $A_{G}$. Let $\varrho=\max\{\rho(G)| G$ be a planar graph of order $n\}$. Our study development on the spectral radius of planar graph is as follows:

\begin{theorem} \label{le6,08}
$\\$
(1) $1.359+\sqrt{2n-\frac{15}{4}}<\varrho< 2+\sqrt{2n-6}$ for $n\geq 10$;
$\\$
(2) $1.472+\sqrt{2n-\frac{15}{4}}<\varrho< 2+\sqrt{2n-6}$ for $n\geq 40$;
$\\$
(3) $1.478+\sqrt{2n-\frac{15}{4}}<\varrho< 2+\sqrt{2n-6}$ for $n\geq 50$.
\end{theorem}

\begin{theorem} \label{le01,01}
Among all planar graphs containing dominating vertex with concrete fixed order $n \geq 48$, $P_{2}\nabla P_{n-2}$ attains the maximum spectral radius.
\end{theorem}

\begin{Oul}
\end{Oul}

The layout of this paper is as follows: section 2 introduces some basic knowledge and working lemmas; section 3 represents our results.

\section{Preliminary}

\ \ \ \ \ For the requirements in the demonstrations afterward, we need some prepares.

From graph theory (see \cite{JBUM} for instance), it is known that:
$\mathrm{(i)}$ a graph is planar graph if and only if it is $K_{3,3}$ and $K_{5}$ minor free; $\mathrm{(ii)}$ a maximal planar graph can be obtained from a non-maximal graph $G$ by adding new edges to $G$;
$\mathrm{(iii)}$ a maximal planar graph $G$ of order $n\geq 3$ is $2$-connected and $\delta\geq 2$;
$\mathrm{(iv)}$ for a planar graph $G$, $\mathbbm{f}$ is invariant, i.e. $\mathbbm{f}$ is constant for different planar embeddings of $G$;
$\mathrm{(v)}$ in any planar embedding of a maximal planar graph $G$ of order $n\geq 3$, $B(f)$ is a 3-cycle (or called triangle) for every face $f$, where $f$ is called triangular face consequently; $\mathrm{(vi)}$ a planar graph is maximal if and only if in its any planar embedding, $B(f)$ is a 3-cycle for every face $f$;
$\mathrm{(vii)}$ for a planar graph $G$ of order $n\geq 3$, $m(G)\leq 3n-6$ with equality if and only if $G$ is maximal; $\mathrm{(viii)}$ for a bipartite planar graph $G$, $m(G)\leq 1$ if $n\leq 2$, and $m(G)\leq 2n-4$ if $n\geq 3$.

For a graph $G$ with vertex set $\{v_{1}$, $v_{2}$, $\ldots$, $v_{n}\}$, a vector $X=(x_{v_1}, x_{v_2}, \ldots, x_{v_n})^T \in R^n$ on $G$ is a vector that
 entry $x_{v_i}$ is mapped to vertex $v_i$ for $i\leq i\leq n$.

Denote by $R^{n}_{++}$ ($R$, $R^{n}_{\geq 0}$) the set of positive real (real, nonnegative real) vectors of dimension $n$. By the famous Perron-Frobenius theorem \cite{H.M, OP}, for $A_{G}$ of a connected graph $G$ of order $n$, we know that there is unique one positive eigenvector $Y=(y_{v_{1}}$, $y_{v_{2}}$, $\ldots$, $y_{v_{n}})^T \in R^{n}_{++}$ corresponding to $\rho(G)$ satisfying $\sum^{n}_{i=1}y^{2}_{v_{i}}= 1$, called $principal$ $eigenvector$; there is unique one positive eigenvector $X=(x_{v_{1}}$, $x_{v_{2}}$, $\ldots$, $x_{v_{n}})^T \in R^{n}_{++}$ corresponding to $\rho(G)$ satisfying $\max\{x_{v_{i}}| 1\leq i\leq n\}= 1$, called $noralized$ $eigenvector$.

Let $A$ be an irreducible nonnegative $n \times n$ real matrix with spectral radius $\rho(A)$ which is the maximum modulus among all eigenvalues of $A$. The following extremal representation (Rayleigh quotient) will be useful:
$$\rho(A)=\max_{X\in R^{n}, X\neq0}\frac{X^{T}AX}{X^{T}X},$$ and if a vector $X$ satisfies that $\frac{X^{T}AX}{X^{T}X}=\rho(A)$, then $AX=\rho(A)X$.

\begin{lemma} \label{le3,01,01}
$\\$
(1) {\bf \cite{JJ,QK}} For a connected graph $G$, $e\notin E(G)$. Then $\rho(G+e)>\rho(G)$.
$\\$
(2) {\bf \cite{JJ}} If $H$ is a proper subgraph of a connected graph $G$, then $\rho(H)<\rho(G)$.
\end{lemma}

\begin{lemma} {\bf\cite{GYHSH}}\label{le5,3,0}
Let $A$ be an irreducible nonnegative square real matrix with order
$n$ and spectral radius $\rho$. If there exists a nonzero
vector $Y=(y_{1}$, $y_{2}$, $\ldots$, $y_{n})^{T}\in R^{n}_{\geq 0}$ and a real coefficient polynomial function $f$ such that $f(A)Y\leq rY$ $(r\in
R)$, then $f(\rho)\leq r$. Similarly, if $f(A)Y\geq rY$ $(r\in
R)$, then $f(\rho)\geq r$.
\end{lemma}

With some modifications of the proof for Lemma \ref{le5,3,0}, we get an improved Lemma \ref{le6,04}.

\begin{lemma} \label{le6,04}
Let $A$ be an irreducible nonnegative square real matrix with order
$n$ and spectral radius $\rho$. $f$ is a real coefficient polynomial function and $Y=(y_{1}$, $y_{2}$, $\ldots$, $y_{n})^{T}\in R^{n}_{\geq 0}$ is a nonzero
vector.

(1) If $f(A)Y\leq rY$ $(r\in
R)$, then $f(\rho)\leq r$ with equality if and only if $f(A)Y= rY$. Moreover, if there is some $y_{i}$ such that $(f(A)Y)_{i}< ry_{i}$ for some $1\leq i\leq n$, then $f(\rho)< r$.

(2) If $f(A)Y\geq rY$ $(r\in
R)$, then $f(\rho)\geq r$ with equality if and only if $f(A)Y= rY$. Moreover, if there is some $y_{i}$ such that $(f(A)Y)_{i}> ry_{i}$ for some $1\leq i\leq n$, then $f(\rho)> r$.
\end{lemma}

\begin{proof}
(1) By Lemma \ref{le5,3,0}, it is enough to consider only the cases that the equalities hold.
We first prove the case for the equality holding in $f(\rho)\leq r$.

Note that $\rho(A^{T})=\rho(A)=\rho$, $A^{T}$ is also irreducible and
nonnegative. Denote by $X=(x_{1}$, $x_{2}$, $\ldots$, $x_{n})^{T}$ the principal eigenvector of $A^{T}$.

Suppose $f(A)Y= rY$. Then $f(\rho)
X^{T}Y=(f(\rho)
X)^{T}Y=(f(A^{T})X)^{T}Y=X^{T}f(A)Y= rX^{T}Y.$ Hence it follows $f(\rho)= r$. Then the sufficiency follows.

Suppose $f(\rho)= r$. If $f(A)Y\neq rY$, then from $f(A)Y\leq rY$, it follows that there is some $y_{i}$ such that $(f(A)Y)_{i}< ry_{i}$.
Then $f(\rho)
X^{T}Y=(f(\rho)
X)^{T}Y=(f(A^{T})X)^{T}Y=X^{T}f(A)Y<r\sum_{j\neq i}x_{j}y_{j}+rx_{i}y_{i}=rX^{T}Y.$
Thus it follows that $f(\rho)< r$, which contradicts $f(\rho)= r$. Then the necessity follows.

Furthermore, from the proof for the necessity, we get that if there is some $y_{i}$ such that $(f(A)Y)_{i}< ry_{i}$ for some $1\leq i\leq n$, then $f(\rho)< r$.

(2) is proved similarly.
This completes the proof. \ \ \ \ \ $\Box$
\end{proof}

\section{Maximum spectral radius of general planar grpah}

\ \ \ \ In this section, we let  $\mathbb{G}_{n}$ be a planar graph of order $n$ satisfying $\rho(\mathbb{G}_{n})=\varrho$; $\mathcal {D}=\{G|\, G$ be a planar graph of order $n\geq 2$ containing dominating vertex$\}$, $\mathcal {G}_{n}\in \mathcal {D}$ satisfy $\rho(\mathcal {G}_{n})=\max\{\rho(G)|\, G\in \mathcal {D}\}$.

Note that a connected planar graph can be obtained from a disconnected planar graph by adding some edges. Using Lemma \ref{le3,01,01} gets the following Lemma \ref{le6,02} immediately.

\begin{lemma} \label{le6,02}
$\\$
(1) $\mathbb{G}_{n}$ is connected.
$\\$
(2) Both $\mathbb{G}_{n}$ and $\mathcal {G}_{n}$ are maximal.
\end{lemma}

\begin{lemma} \label{le6,02,01}
$\\$
(1) In any embedding of $\mathbb{G}_{n}$, for a vertex $v_{i}$, then in $\mathbb{G}^{v_{i}}_{n}$, there is a cycle formed by all the vertices in $N_{\mathbb{G}}(v_{i})$, denoted by $\mathbb{C}_{i}=v_{i_{1}}v_{i_{2}}\cdots v_{i_{j}}v_{i_{1}}$ where $j=d_{\mathbb{G}}(v_{i})$, $N_{\mathbb{G}}(v_{i})=\{v_{i_{1}}$, $v_{i_{2}}$, $\ldots$, $v_{i_{j}}\}$, and $v_{i}\bigtriangledown \mathbb{C}_{i}$ is a subgraph.
$\\$
(2) In any embedding of $\mathcal {G}_{n}$, for a vertex $v_{i}$, then in $\mathcal {G}^{v_{i}}_{n}$, there is a cycle formed by all the vertices in $N_{\mathcal {G}_{n}}(v_{i})$, denoted by $\mathbb{C}_{i}=v_{i_{1}}v_{i_{2}}\cdots v_{i_{j}}v_{i_{1}}$ where $j=d_{\mathcal {G}_{n}}(v_{i})$, $N_{\mathcal {G}_{n}}(v_{i})=\{v_{i_{1}}$, $v_{i_{2}}$, $\ldots$, $v_{i_{j}}\}$, and $v_{i}\bigtriangledown \mathbb{C}_{i}$ is a subgraph.
\end{lemma}

\begin{proof}
Suppose $\widetilde{\mathbb{G}_{n}}$ is a planar embedding of $\mathbb{G}_{n}$. For convenience, we use some new labels for the vertices of $\mathbb{G}_{n}$. In $\widetilde{\mathbb{G}_{n}}$, along clockwise direction, denote by $w_{i_{1}}$, $w_{i_{2}}$, $\ldots$, $w_{i_{j}}$ the neighbors of $v_{i}$. Note that in $\widetilde{\mathbb{G}_{n}}$, every face of $\mathbb{G}_{n}$ is a triangular face; every pair of adjacent edges $v_{i}w_{i_{\iota}}$ and $v_{i}w_{i_{\iota+1}}$ are in a common face (a triangular face) for each $1\leq \iota \leq j-1$; $v_{i}w_{i_{j}}$ and $v_{i}w_{i_{1}}$ are in a common face. Thus in $\widetilde{\mathbb{G}_{n}}$, $w_{i_{1}}w_{i_{2}}\cdots w_{i_{j}}w_{i_{1}}$ forms a cycle, denoted by $\mathbb{C}_{i}$, and $v_{i}\bigtriangledown \mathbb{C}_{i}$ is a subgraph. Let $v_{i_{\iota}}=w_{i_{\iota}}$ for $1\leq \iota \leq j$. Then (1) follows. (2) is proved similarly. This completes the proof. \ \ \ \ \ $\Box$
\end{proof}

Let $r_{i}(A)$ be the ith row sum of a matrix $A$, $r_{v}(A_{G})$ be the row sum of $A_{G}$ corresponding to vertex $v$ for a graph $G$.

\begin{lemma} \label{le6,17}
 Let $G$ be a planar graph on $n\geq 7$ vertices. Then $\rho(G)\leq 2+\sqrt{2n-6}$.
\end{lemma}

\begin{proof}
For any vertex $v$ in $G$, noting that $G^{v}$ is planar, we have $e(N_{G}[v]) \leq 3(d_{v}+1) - 6$, implying that $e(N_{G}(v)) \leq
 2d_{v} - 3$. Note that $G(N_{G}(v), V(G)\setminus N_{G}(v))$ is bipartite. If $2 \leq d_{v} \leq n-2$, then
$e(N_{G}(v), V(G)\setminus N_{G}(v)) \leq 2n -4$. Otherwise, if $d_{v} =1$ or $d_{v} =n-1$, then $e(N_{G}(v), V(G)\setminus N_{G}(v)) \leq n-1\leq 2n -4$ for $n\geq 3$. Thus
$$r_{v}(A^{2}_{G})=\sum_{u\sim v}d_{u}=2e(N_{G}(v))+e(N_{G}(v), V(G)\setminus N_{G}(v))\leq 4d_{v} - 6+2n -4=2n+4d_{v}-10.$$

It follows that $r_{v}(A^{2}_{G})\leq 2n+4d_{v}-10$, i.e. $r_{v}(A^{2}_{G})\leq 2n+4r_{v}(A_{G})-10$. Thus for vector $Y=(1$, $1$, $\ldots$, $1)^{T}$, we have $(A^{2}_{G}-4A_{G}-2n+10)Y\leq 0Y$. Then by Lemma \ref{le6,04},  it follows that $\rho^{2}(G)-4\rho(G)-2n+10\leq 0$, implying that $\rho(G)\leq 2+\sqrt{2n-6}$.
This completes the proof. \ \ \ \ \ $\Box$
\end{proof}

\setlength{\unitlength}{0.65pt}
\begin{center}
\begin{picture}(632,313)
\put(170,284){\circle*{4}}
\put(170,118){\circle*{4}}
\qbezier(170,284)(170,201)(170,118)
\put(86,86){\circle*{4}}
\put(171,244){\circle*{4}}
\put(170,156){\circle*{4}}
\put(254,86){\circle*{4}}
\qbezier(86,86)(170,86)(254,86)
\qbezier(170,284)(84,267)(86,86)
\qbezier(170,284)(265,268)(254,86)
\qbezier(171,244)(128,165)(86,86)
\qbezier(171,244)(212,165)(254,86)
\qbezier(170,156)(128,121)(86,86)
\qbezier(170,156)(212,121)(254,86)
\qbezier(170,118)(128,102)(86,86)
\qbezier(170,118)(212,102)(254,86)
\qbezier(170,118)(67,0)(38,81)
\qbezier(170,284)(0,313)(38,81)
\put(518,288){\circle*{4}}
\put(518,111){\circle*{4}}
\qbezier(518,288)(518,200)(518,111)
\put(435,73){\circle*{4}}
\put(601,73){\circle*{4}}
\qbezier(435,73)(518,73)(601,73)
\put(518,252){\circle*{4}}
\put(518,145){\circle*{4}}
\qbezier(518,288)(419,251)(435,73)
\qbezier(518,288)(612,255)(601,73)
\qbezier(518,252)(476,163)(435,73)
\qbezier(518,252)(559,163)(601,73)
\qbezier(518,145)(476,109)(435,73)
\qbezier(518,145)(559,109)(601,73)
\qbezier(518,111)(476,92)(435,73)
\qbezier(518,111)(559,92)(601,73)
\put(80,74){$v_{1}$}
\put(247,74){$v_{2}$}
\put(166,104){$v_{3}$}
\put(174,156){$v_{4}$}
\put(175,244){$v_{n-1}$}
\put(168,289){$v_{n}$}
\put(156,38){$\Phi$}
\put(424,63){$v_{1}$}
\put(603,63){$v_{2}$}
\put(513,98){$v_{3}$}
\put(522,147){$v_{4}$}
\put(522,250){$v_{n-1}$}
\put(514,294){$v_{n}$}
\put(513,38){$\mathcal {H}$}
\put(293,4){Fig. 3.1. $\Phi,\ \mathcal {H}$}
\end{picture}
\end{center}

Let $n\geq 4$, $P=v_{1}v_{2}$, $C=v_{3}v_{4}\cdots v_{n}v_{3}$, $\Phi=P\bigtriangledown C$, $\mathcal {H}=\Phi-v_{3}v_{n}$ (see Fig. 3.1).

\begin{lemma} \label{le6,06}
$\rho(\Phi)=\frac{3}{2}+\sqrt{2n-\frac{15}{4}}$ for $n\geq 5$; $\rho(\mathcal {H})\geq \frac{3}{2}+\sqrt{2n-\frac{15}{4}}-\frac{2}{2n-\frac{15}{4}-\frac{\sqrt{2n-\frac{15}{4}}}{2}}$ for $n\geq 5$.
\end{lemma}

\begin{proof}
Let $Y=(y_{v_{1}}$, $y_{v_{2}}$, $\ldots$, $y_{v_{n}})^T \in R^{n}_{++}$ be the principal eigenvector of $\Phi$. By symmetry, it follows that $y_{v_{1}}=y_{v_{2}}$, $y_{v_{3}}=y_{v_{4}}=\cdots=y_{v_{n}}$. Note that $$\left \{\begin{array}{ll}
 \rho(\Phi)y_{v_{1}}=y_{v_{2}}+(n-2)y_{v_{3}}=y_{v_{1}}+(n-2)y_{v_{3}},\\
\\ \rho(\Phi)y_{v_{3}}=y_{v_{1}}+y_{v_{2}}+y_{v_{4}}+y_{v_{n}}=2y_{v_{1}}+2y_{v_{3}},\end{array}\right.$$
i.e. $$\left \{\begin{array}{ll}
 (\rho(\Phi)-1)y_{v_{1}}=(n-2)y_{v_{3}},\\
\\ (\rho(\Phi)-2)y_{v_{3}}=2y_{v_{1}}.\end{array}\right.$$
It follows that $\rho^{2}(\Phi)-3\rho(\Phi)-2n+6=0$ and $\rho(\Phi)=\frac{3}{2}+\sqrt{2n-\frac{15}{4}}$.

From $2y_{v_{1}}=(\rho(\Phi)-2)y_{v_{3}}=(\sqrt{2n-\frac{15}{4}}-\frac{1}{2})y_{v_{3}}$ and $\sum^{n}_{i=1}y^{2}_{v_{i}}=1$, it follows that
$y^{2}_{v_{3}}=\frac{1}{2n-\frac{15}{4}-\frac{\sqrt{2n-\frac{15}{4}}}{2}}$.

Then we get $$\rho(\mathcal {H})\geq Y^{T}A_{\mathcal {H}}Y=Y^{T}A_{\Phi}Y-2y_{v_{3}}y_{v_{n}}=\frac{3}{2}+\sqrt{2n-\frac{15}{4}}-2y^{2}_{v_{3}}
=\frac{3}{2}+\sqrt{2n-\frac{15}{4}}-\frac{2}{2n-\frac{15}{4}-\frac{\sqrt{2n-\frac{15}{4}}}{2}}.$$
This completes the proof. \ \ \ \ \ $\Box$
\end{proof}

Note that if positive integer $n\geq 2$, then $f(n)=\frac{3}{2}+\sqrt{2n-\frac{15}{4}}-\frac{2}{2n-\frac{15}{4}-\frac{\sqrt{2n-\frac{15}{4}}}{2}}$ and $2n-\frac{15}{4}-\frac{\sqrt{2n-\frac{15}{4}}}{2}$ strictly increase with respect to $n\geq 4$ increasing. Thus we get the following Corollary \ref{le6,06,01}.

\begin{corollary} \label{le6,06,01}
$\\$
(1) $\rho(\mathcal {H})> 1.359+\sqrt{2n-\frac{15}{4}}$ for $n\geq 10$;
$\\$
(2) $\rho(\mathcal {H})> 1.472+\sqrt{2n-\frac{15}{4}}$ for $n\geq 40$;
$\\$
(3) $\rho(\mathcal {H})> 1.478+\sqrt{2n-\frac{15}{4}}$ for $n\geq 50$.
\end{corollary}

Note the maximality of $\rho(\mathbb{G}_{n})$ and $\mathcal {H}$ is a planar graph. Combining Lemma \ref{le6,17} and Corollary \ref{le6,06,01}, we get the following Theorem \ref{le6,18}.

\begin{theorem} \label{le6,18}
$\\$
(1) $1.359+\sqrt{2n-\frac{15}{4}}<\rho(\mathbb{G}_{n})< 2+\sqrt{2n-6}$ for $n\geq 10$;
$\\$
(2) $1.472+\sqrt{2n-\frac{15}{4}}<\rho(\mathbb{G}_{n})< 2+\sqrt{2n-6}$ for $n\geq 40$;
$\\$
(3) $1.478+\sqrt{2n-\frac{15}{4}}<\rho(\mathbb{G}_{n})< 2+\sqrt{2n-6}$ for $n\geq 50$.
\end{theorem}

\begin{Proof}
This theorem follows from \ref{le6,18}.
This completes the proof. \ \ \ \ \ $\Box$
\end{Proof}

Note that $\mathcal {H}$ is also a planar graph containing dominating vertex. The next Corollary \ref{le6,07} follows from Corollary \ref{le6,06,01} immediately.

\begin{corollary} \label{le6,07}
$\\$
(1) $\rho(\mathcal {G}_{n})> 1.359+\sqrt{2n-\frac{15}{4}}$ for $n\geq 10$;
 $\\$
(2) $\rho(\mathcal {G}_{n})> 1.472+\sqrt{2n-\frac{15}{4}}$ for $n\geq 40$;
$\\$
(3) $\rho(\mathcal {G}_{n})> 1.478+\sqrt{2n-\frac{15}{4}}$ for $n\geq 50$.
\end{corollary}

Hereafter, we suppose that $V(\mathcal {G}_{n})=\{v_{1}$, $v_{2}$, $\ldots$, $v_{n}\}$, and $deg_{\mathcal {G}_{n}}(v_{1})=n-1$. Let $X=(x_{v_{1}}$, $x_{v_{2}}$, $\ldots$, $x_{v_{n}})^T \in R^{n}_{++}$ be the normalized eigenvector of $\mathcal {G}_{n}$.

\begin{lemma} \label{le6,07,0}
$x_{v_{1}}=1$ if $n\geq 3$.
\end{lemma}

\begin{proof}
For $v_{1}$, and any vertex $v_{i}$ where $2\leq i\leq n$, using $\rho(\mathcal {G}_{n})(x_{v_{1}}-x_{v_{i}})=x_{v_{i}}-x_{v_{1}}+\sum_{3\leq j\leq n}x_{v_{j}}-\sum_{v_{j}\sim v_{i}}x_{v_{j}}$ gets that $(\rho(\mathcal {G}_{n})-1)(x_{v_{1}}-x_{v_{i}})=\sum_{3\leq j\leq n}x_{v_{j}}-\sum_{v_{j}\sim v_{i}}x_{v_{j}}\geq 0$. Note that $\mathcal {G}_{n}$ is maximal and $\mathcal {G}_{n}$ contains at least one cycle in $\mathcal {G}_{n}$. Then $\rho(\mathcal {G}_{n})\geq 2$. Thus it follows that $x_{v_{1}}-x_{v_{i}}=\frac{\sum_{3\leq j\leq n}x_{v_{j}}-\sum_{v_{j}\sim v_{i}}x_{v_{j}}}{\rho(\mathcal {G}_{n})-1}\geq 0$. As a result, it follows that $x_{v_{1}}\geq x_{v_{i}}$. Note that $X$ is the normalized eigenvector of $\mathcal {G}_{n}$ and the arbitrariness of $v_{i}$. Then we get that $x_{v_{1}}=1$.
This completes the proof. \ \ \ \ \ $\Box$
\end{proof}

Next, we let $1=x_{v_{1}}\geq x_{v_{2}}\geq \cdots \geq x_{v_{n}}$ for convenience.

\begin{lemma} \label{le6,14}
$x_{v_{2}}> \frac{1}{4}$ if $n\geq 10$.
\end{lemma}

\begin{proof}
Let $H=\mathcal {G}_{n}-v_{1}$. Note that $\mathcal {G}_{n}$ is maximal, $x_{v_{1}}=1$ and note that for $n\geq 10$, $e(\mathcal {G}_{n})=3n-6$, $d_{\mathcal {G}_{n}}(v_{1})=n-1$, $\rho(\mathcal {G}_{n})> 1.359+\sqrt{2n-\frac{15}{4}}$. Then $\sum_{i\geq 2}d_{H}(v_{i})\leq 2(2n-5)$, and
$$\rho^{2}(\mathcal {G}_{n})x_{v_{1}}=\sum_{v_{j}\sim v_{1}} \rho(\mathcal {G}_{n}) x_{v_{j}}=(n-1)x_{v_{1}}+\sum_{i\geq 2}d_{H}(v_{i})x_{v_{i}}\leq (n-1)x_{v_{1}}+2(2n-5)x_{v_{2}}$$
$$\Longrightarrow (\rho^{2}(\mathcal {G}_{n})-n+1)x_{v_{1}}\leq (4n-10)x_{v_{2}}\hspace{3.75cm}$$ $$\Longrightarrow \frac{\rho^{2}(\mathcal {G}_{n})-n+1}{4n-10}x_{v_{1}}\leq x_{v_{2}}\hspace{5.35cm}$$ $$\hspace{1.5cm}\Longrightarrow \frac{1.846+2.718\sqrt{2n-\frac{15}{4}}+n-\frac{11}{4}}{4n-10}x_{v_{1}}< x_{v_{2}}\hspace{3.4cm} (\ast 1)$$
$$\Longrightarrow x_{v_{2}}> \frac{1}{4}\hspace{7.8cm}$$

This completes the proof. \ \ \ \ \ $\Box$
\end{proof}

Let $f_{2}(n)=\frac{1.846+2.718\sqrt{2n-\frac{15}{4}}+n-\frac{11}{4}}{4n-10}$. Taking derivation with respect to $n$ gets $$f'_{2}(n)=\frac{(\frac{2.718}{2\sqrt{2n-\frac{15}{4}}}+1)(4n-10)-4(1.846+2.718\sqrt{2n-\frac{15}{4}}+n-\frac{11}{4})}{(4n-10)^{2}}< 0.$$ It follows that $f_{2}(n)=\frac{1.846+2.718\sqrt{2n-\frac{15}{4}}+n-\frac{11}{4}}{4n-10}$ decreases strictly with respect to $n\geq 1$ increasing. Then using inequality $(\ast 1)$ gets the following Lema \ref{le6,14,01}.

\begin{lemma} \label{le6,14,01}
$\\$
(1) $x_{v_{2}}> 0.3987$ if $n\leq 50$;
$\\$
(2) $x_{v_{2}}> 0.375$ if $n\leq 80$.
\end{lemma}

\setlength{\unitlength}{0.7pt}
\begin{center}
\begin{picture}(300,126)
\put(24,58){\circle*{4}}
\put(102,95){\circle*{4}}
\qbezier(24,58)(63,77)(102,95)
\put(225,107){\circle*{4}}
\qbezier(102,95)(163,101)(225,107)
\put(295,48){\circle*{4}}
\qbezier(225,107)(260,78)(295,48)
\put(170,31){\circle*{4}}
\qbezier(102,95)(136,63)(170,31)
\qbezier(295,48)(232,40)(170,31)
\qbezier(24,58)(97,45)(170,31)
\qbezier(102,95)(198,72)(295,48)
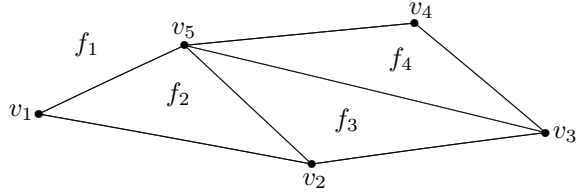
\put(-30,-9){Fig. 3.2. an outerplanar graph $G$ and its an OP-embedding $\widehat{G}$}
\put(8,57){$v_{1}$}
\put(165,20){$v_{2}$}
\put(299,45){$v_{3}$}
\put(222,112){$v_{4}$}
\put(95,101){$v_{5}$}
\put(43,92){$f_{1}$}
\put(92,65){$f_{2}$}
\put(182,51){$f_{3}$}
\put(211,84){$f_{4}$}
\end{picture}
\end{center}

A simple graph $G$ is $outerplanar$ if it has an embedding in
the plane, called $outerplane$-$embedding$ (written as $OP$-$embedding$ for short hereafter), denoted by $\widehat{G}$, so that every vertex lies on the boundary of the unbounded (outer) face. Similar to a planar embedding of a planar graph, an OP-embedding of an outerplanar graph $G$ partitions the plane into a number of edgewise-connected
$faces$. Among the faces of $\widehat{G}$ for an outerplanar graph $G$, the outer one is called the $outer$ $face$ (see $f_{1}$ in Fig. 3.2 for example), and any one of other faces is called the $inner$ $face$ (see $f_{2}$, $f_{3}$, $f_{4}$ in Fig. 3.2 for example). Similar to the notations for a planar graph, we denote by $B(f)$ the boundary of a face $f$ in $\widehat{G}$ of an outerplanar graph $G$. It can be seen that $B(f)$ of a face $f$ is a circuit. For an outerplanar graph $G$, we denote by $O_{\widehat{G}}$ the outer face of its an OP-embedding $\widehat{G}$. As shown in Fig. 3.2, we can see that $O_{\widehat{G}}=f_{1}$, $B(f_{1})=v_{1}v_{2}v_{3}v_{4}v_{5}v_{1}$.

A simple outerplanar graph is (edge) $maximal$ if no edge can be added to
the graph without violating its simplicity, or outerplanarity.

Recall some facts (known results in graph theory, see \cite{JBUM} for instance) about outerplanar graph that $\mathrm{(i)}$ a maximal outerplanar graph can be obtained from a non-maximal outerplanar graph $G$ by adding new edges to $G$; $\mathrm{(ii)}$ in a maximal outerplanar graph $G$ of order $n\geq 3$, $\delta(G)= 2$ and there are at least two vertices with degree $2$; $\mathrm{(iii)}$ in an OP-embedding $\widehat{G}$ of a maximal outerplanar graph $G$, the boundary of every inner face is a 3-cycle (or called triangle), where every inner face is called a triangular face similarly, the $B(O_{\widehat{G}})$ is a Hamilton cycle of $G$; $\mathrm{(iv)}$ an outerplanar graph is maximal if and only if there is an OP-embedding $\widehat{G}$ such that the boundary of every inner face is a 3-cycle and the $B(O_{\widehat{G}})$ is a Hamilton cycle of $G$; $\mathrm{(v)}$ for an outerplanar graph $G$ of order $n$, $e(G)\leq 2n-3$ with equality if and only if $G$ is maximal.

From Lemma \ref{le6,02,01} and definition of outerplanar graph, it follows that both $\mathcal {G}_{n}-v_{1}$ and $\mathcal {G}_{n}-\{v_{1}, v_{2}\}$ are outerplanar.

Let $\mathbb{L} =\{v_{1}, v_{2}\}$, $\mathbb{S}=V(\mathcal {G}_{n})\setminus\mathbb{L}$.

\begin{lemma} \label{le6,12}
If $n\geq 81$, then for any $i\geq 3$, we have $x_{v_{i}}<\frac{x_{v_{2}}}{2}$.
\end{lemma}

\begin{proof}
We prove this result by contradiction.
Suppose $x_{v_{i}}\geq\frac{x_{v_{2}}}{2}$ for $i\geq 3$. Then $$(n-2)\frac{x_{v_{2}}}{2}\leq\sum_{v_{z}\in \mathbb{S}}x_{v_{z}}=\frac{1}{\rho(\mathcal {G}_{n})}\sum_{v_{z}\in \mathbb{S}}\rho(\mathcal {G}_{n}) x_{v_{z}}\hspace{7.8cm}$$$$\leq \frac{1}{1.478+\sqrt{2n-\frac{15}{4}}}\sum_{v_{z}\in \mathbb{S}}\rho(\mathcal {G}_{n}) x_{v_{z}}\hspace{2.5cm}$$$$=\frac{1}{1.478+\sqrt{2n-\frac{15}{4}}}\sum_{v_{z}\in \mathbb{S}}\sum_{v_{y}\sim v_{z}} x_{v_{y}}\hspace{2.5cm}$$ $$\hspace{0.9cm}= \frac{1}{1.478+\sqrt{2n-\frac{15}{4}}}\sum_{v_{z}\in \mathbb{S}}(\sum_{v_{y}\sim v_{z}, v_{y}\in \mathbb{S}} x_{v_{y}}+\sum_{v_{y}\sim v_{z}, v_{y}\in \mathbb{L}} x_{v_{y}})$$
$$\hspace{3.45cm}\leq \frac{1}{1.478+\sqrt{2n-\frac{15}{4}}}(x_{v_{2}}2e(\mathbb{S})+(n-2)x_{v_{1}}+(n-2)x_{v_{2}})\hspace{1.9cm}(\ast 2)$$
Note that $\mathcal {G}_{n}[\mathbb{S}]=\mathcal {G}_{n}-\{v_{1}, v_{2}\}$\ is\ outerplanar. Then
$$(\ast 2)\leq \frac{1}{1.478+\sqrt{2n-\frac{15}{4}}}(x_{v_{2}}(2(n-2)-3)+(n-2)x_{v_{1}}+(n-2)x_{v_{2}})\hspace{0.15cm}$$
$$=\frac{1}{1.478+\sqrt{2n-\frac{15}{4}}}(x_{v_{2}}(3n-9)+n-2)\hspace{3.15cm}$$
and $$\hspace{0.3cm}(n-2)\frac{x_{v_{2}}}{2}<\frac{1}{1.478+\sqrt{2n-\frac{15}{4}}}(x_{v_{2}}(3n-9)+n-2).\hspace{5.3cm}(\ast 3)$$
Note that if $n\geq 81$, then $1.478+\sqrt{2n-\frac{15}{4}}> 14.05$. Combining Lemma \ref{le6,14}, from ($\ast 3$), it follows that $$\frac{x_{v_{2}}}{2}\leq\frac{\frac{1}{1.478+\sqrt{2n-\frac{15}{4}}}(x_{v_{2}}(3n-9)+n-2)}{n-2}<
\frac{3x_{v_{2}}+1}{1.478+\sqrt{2n-\frac{15}{4}}}<\frac{x_{v_{2}}}{2}.\hspace{2.8cm}(\ast 4)$$
It is a contradiction. Thus our result follows.
This completes the proof. \ \ \ \ \ $\Box$
\end{proof}

Combining Lemma \ref{le6,14,01} and inequality $(\ast 4)$, and repeating the process in proof of Lemma \ref{le6,12}, we get the following Lemma \ref{le6,15}.

\begin{lemma} \label{le6,15}
If $51\leq n\leq 80$, then for any $i\geq 3$, we have $x_{v_{i}}<\frac{x_{v_{2}}}{2}$.
\end{lemma}

Combining $\rho(\mathcal {G}_{n})> 1.472+\sqrt{2n-\frac{15}{4}}$ for $n\geq 40$ in Lemma \ref{le6,07} and inequality $(\ast 4)$, and repeating the process in proof of Lemma \ref{le6,12} by replacing $1.478+\sqrt{2n-\frac{15}{4}}$ with $1.472+\sqrt{2n-\frac{15}{4}}$, we get the following Lemma \ref{le6,15}.

\begin{lemma} \label{le6,15}
If $48\leq n\leq 50$, then for any $i\geq 3$, we have $x_{v_{i}}<\frac{x_{v_{2}}}{2}$.
\end{lemma}

\begin{lemma} \label{le6,13}
If $n\geq 48$, then $d_{\mathcal {G}_{n}}(v_{2})=n-1$.
\end{lemma}

\begin{proof}
We prove this result by contradiction. Suppose $d_{\mathcal {G}_{n}}(v_{2})\leq n-2$, $\widetilde{\mathcal {G}_{n}}$ is a planar embedding of $\mathcal {G}_{n}$. For convenience, we use some new labels for the vertices of $\mathcal {G}_{n}$. By Lemma \ref{le6,02,01}, in $\widetilde{\mathcal {G}_{n}}$, suppose $\mathbb{C}_{1}=v_{2}w_{1}w_{2}\cdots w_{n-2}v_{2}$ is a cycle where $N_{\mathcal {G}_{n}}(v_{1})=\{v_{2}$, $w_{1}$, $w_{2}$, $\ldots$, $w_{n-2}\}$, and $v_{1}\bigtriangledown \mathbb{C}_{1}$ is a subgraph. Suppose $v_{2}$, $w_{1}$, $w_{2}$, $\ldots$, $w_{n-2}$ are distributed along clockwise direction around vertex $v_{1}$ in $\widetilde{\mathcal {G}_{n}}$. Suppose $N_{\mathcal {G}_{n}}(v_{2})=\{v_{1}$, $w_{i_{1}}$, $w_{i_{2}}$, $\ldots$, $w_{i_{\eta}}\}$ and $1=i_{1}<i_{2}<\cdots<i_{\eta}=n-2$. For convenience to distinguish $i_{t}+1$ and $i_{t+1}$ ($i_{t}-1$ and $i_{t-1}$) easily, we use $w_{i_{(t+1)}}$ for $w_{i_{t+1}}$ ($w_{i_{(t-1)}}$ for $w_{i_{t-1}}$) hereafter. Because $d_{\mathcal {G}_{n}}(v_{2})\leq n-2$, there exist $w_{i_{j}}$ and $w_{i_{(j+1)}}$ such that
$i_{(j+1)}-i_{j}\geq 2$ where $1\leq j\leq \eta-1$. This means that none of vertices $w_{i_{j}+1}$, $w_{i_{j}+2}$, $\ldots$, $w_{i_{(j+1)}-1}$ is adjacent to $v_{2}$.

Note that $\mathcal {G}_{n}$ is maximal, and in $\widetilde{\mathcal {G}_{n}}$, $v_{2}w_{i_{j}}$, $v_{2}w_{i_{(j+1)}}$ are in a common triangular face. Then $w_{i_{j}}w_{i_{(j+1)}}\in E(\mathcal {G}_{n})$.
By the facts about maximal outerplanar graph, $H=\mathcal {G}_{n}[\{w_{i_{j}}$, $w_{i_{j}+1}$, $w_{i_{j}+2}$, $\ldots$, $w_{i_{(j+1)}-1}$, $w_{i_{(j+1)}}\}]$ is a maximal outerplanar graph, $C=w_{i_{j}}w_{i_{j}+1}w_{i_{j}+2}\cdots w_{i_{(j+1)}-1}w_{i_{(j+1)}}w_{i_{j}}$ is the Hamilton cycle of $H$.

{\bf Claim} If $i_{(j+1)}-i_{j}\geq 3$, then at least one of $d_{H}(w_{i_{j}})$, $d_{H}(w_{i_{(j+1)}})$ is more than $2$. With respect to $\widetilde{\mathcal {G}_{n}}$, noting the maximality $\mathcal {G}_{n}$ and $H$, we know that every inner face of $H$ is a triangular face. Thus $w_{i_{j}}w_{i_{(j+1)}}$ is in a inner triangular face of $H$. Then in $H$, there is a vertex $w_{\varsigma}$ together with $w_{i_{j}}$, $w_{i_{(j+1)}}$ forming a $3$-cycle $w_{i_{j}}w_{\varsigma}w_{i_{(j+1)}}w_{i_{j}}$, where $i_{j}+1\leq \varsigma\leq i_{(j+1)}-1$. Because $i_{(j+1)}-i_{j}\geq 3$, then at least one of $\varsigma-i_{j}\geq 2$, $i_{(j+1)}-\varsigma\geq 2$ holds. Without loss of generality, suppose $\varsigma-i_{j}\geq 2$. Note that $\{w_{i_{j}}w_{i_{(j+1)}}$, $w_{i_{j}}w_{i_{j}+1}$, $w_{i_{j}}w_{\varsigma}\} \subseteq E(\mathcal {G}_{n})$. Thus $d_{H}(w_{i_{j}})\geq 3$. Consequently, the claim holds.

\setlength{\unitlength}{0.65pt}
\begin{center}
\begin{picture}(662,591)
\put(359,364){\circle*{4}}
\put(85,416){\circle*{4}}
\put(628,416){\circle*{4}}
\qbezier(85,416)(465,691)(628,416)
\put(196,416){\circle*{4}}
\qbezier(85,416)(177,454)(196,416)
\put(316,416){\circle*{4}}
\qbezier(85,416)(228,488)(316,416)
\put(398,416){\circle*{4}}
\qbezier(85,416)(385,558)(398,416)
\put(546,416){\circle*{4}}
\qbezier(85,416)(433,612)(546,416)
\qbezier(85,416)(356,416)(628,416)
\put(359,416){\circle*{4}}
\put(140,416){\circle*{4}}
\qbezier(85,416)(87,361)(359,364)
\qbezier(359,364)(602,358)(628,416)
\put(353,351){$v_{1}$}
\put(67,416){$v_{2}$}
\put(125,405){$w_{1}$}
\put(197,422){$w_{i_{2}}$}
\put(298,405){$w_{i_{j}}$}
\put(346,425){$w_{i_{j}+1}$}
\put(400,405){$w_{i_{(j+1)}}$}
\put(545,404){$w_{i_{(\eta-1)}}$}
\put(633,413){$w_{n-2}$}
\qbezier(316,416)(361,485)(398,416)
\qbezier(140,416)(160,378)(359,364)
\qbezier(196,416)(199,389)(359,364)
\put(397,414){\circle*{4}}
\qbezier(359,364)(401,379)(397,414)
\put(351,319){$\mathcal {G}_{n}$}
\put(86,105){\circle*{4}}
\put(194,105){\circle*{4}}
\qbezier(86,105)(171,153)(194,105)
\put(318,105){\circle*{4}}
\qbezier(86,105)(256,205)(318,105)
\put(408,105){\circle*{4}}
\qbezier(86,105)(374,271)(408,105)
\qbezier(318,105)(364,176)(408,105)
\put(546,105){\circle*{4}}
\qbezier(86,105)(424,330)(546,105)
\put(626,105){\circle*{4}}
\qbezier(86,105)(478,396)(626,105)
\qbezier(316,416)(337,390)(359,364)
\qbezier(86,105)(202,105)(318,105)
\qbezier(408,105)(517,105)(626,105)
\qbezier(86,105)(93,49)(357,48)
\qbezier(357,48)(586,38)(626,105)
\put(351,34){$v_{1}$}
\put(71,110){$v_{2}$}
\put(140,105){\circle*{4}}
\qbezier(140,105)(143,70)(357,48)
\qbezier(194,105)(200,73)(357,48)
\qbezier(318,105)(337,77)(357,48)
\qbezier(408,105)(403,61)(357,48)
\qbezier(359,364)(526,369)(546,416)
\put(123,95){$w_{1}$}
\put(195,111){$w_{i_{2}}$}
\put(304,92){$w_{i_{j}}$}
\put(409,93){$w_{i_{(j+1)}}$}
\put(547,94){$w_{i_{(\eta-1)}}$}
\put(631,103){$w_{n-2}$}

\put(36,105){\circle*{4}}
\put(22,114){$w_{i_{j}+1}$}
\put(357,48){\circle*{4}}
\qbezier(357,48)(41,31)(36,105)
\qbezier(359,416)(359,390)(359,364)
\put(357,49){\circle*{4}}
\put(546,105){\circle*{4}}
\qbezier(357,49)(515,53)(546,105)
\qbezier(36,105)(61,105)(86,105)
\put(351,10){$\mathcal {G}'_{n}$}
\put(305,-17){Fig. 3.3. $\mathcal {G}_{n}$, $\mathcal {G}'_{n}$}
\end{picture}
\end{center}

{\bf Case 1} $i_{(j+1)}-i_{j}= 2$. Note that $H$ is maximal outerplanar, and note the facts about maximal outerplanar graph. It follows that $w_{i_{j}}w_{i_{j}+1}w_{i_{(j+1)}}w_{i_{j}}$ is a $3$-cycle. It means that $d_{H}(w_{i_{j}+1})=2$. Let $\mathcal {G}'_{n}=\mathcal {G}_{n}-w_{i_{j}}w_{i_{j}+1}-w_{i_{j}+1}w_{i_{(j+1)}}+v_{2}w_{i_{j}+1}$ (see Fig. 3.3). Then using Lemmas \ref{le6,12}-\ref{le6,15} gets $X^{T}A_{\mathcal {G}^{'}_{n}}X-X^{T}A_{\mathcal {G}_{n}}X=x_{w_{i_{j}+1}}(x_{v_{2}}-x_{w_{i_{j}}}-x_{w_{i_{(j+1)}}})>0$. It follows that $\rho(\mathcal {G}'_{n})> \rho(\mathcal {G}_{n})$, which contradicts the maximality of $\rho(\mathcal {G}_{n})$.

{\bf Case 2} $i_{(j+1)}-i_{j}\geq 3$. By our above Claim, and the fact that a maximal outerplanar graph of order at least $3$ has at least 2 vertices with degree $2$, there is a vertex $w_{z}$ with $d_{H}(w_{z})=2$ where $i_{j}+1\leq z\leq i_{(j+1)}-1$. Let $\mathcal {G}'_{n}=\mathcal {G}_{n}-w_{z}w_{z+1}-w_{z}w_{z-1}+v_{2}w_{z}$. Similar to Case 1, we get a contradiction that $\rho(\mathcal {G}'_{n})> \rho(\mathcal {G}_{n})$ contradicts the maximality of $\rho(\mathcal {G}_{n})$.

The above narrations demonstrate that $d_{\mathcal {G}_{n}}(v_{2})=n-1$. This completes the proof. \ \ \ \ \ $\Box$
\end{proof}

\begin{theorem} \label{le6,16}
If $n\geq 48$, then $\mathcal {G}_{n}\cong P_{2}\bigtriangledown P_{n-2}$.
\end{theorem}

\begin{proof}
By Lemma \ref{le6,02,01}, in a planar embedding $\widetilde{\mathcal {G}_{n}}$ of $\mathcal {G}_{n}$, suppose $\mathbb{C}_{1}=v_{2}w_{1}w_{2}\cdots w_{\kappa}v_{2}$ is a cycle where $N_{\mathcal {G}_{n}}(v_{1})=\{v_{2}$, $w_{1}$, $w_{2}$, $\ldots$, $w_{n-2}\}$, and $v_{1}\bigtriangledown \mathbb{C}_{1}$ is a subgraph. Let $\mathbb{P}_{2}=v_{1}v_{2}$, $\mathbb{P}_{n-2}=w_{1}w_{2}\cdots w_{n-2}$. By Lemma \ref{le6,13}, it follows that $\mathcal {G}_{n}=\mathbb{P}_{2}\bigtriangledown \mathbb{P}_{n-2}$. Thus the result follows as desired.
This completes the proof. \ \ \ \ \ $\Box$
\end{proof}

\begin{prooff}
This theorem follows from Theorem \ref{le6,16}. This completes the proof. \ \ \ \ \ $\Box$
\end{prooff}

\small {

}

\end{document}